# Approach to design of Nonlinear Robust Control in a Class of Structurally Stable Functions

Viktor Ten

**Abstract**

An approach to stabilization of control systems with ultimately wide ranges of uncertainly disturbed parameters is offered. The method relies on using of nonlinear structurally stable functions from catastrophe theory as controllers. Analytical part presents an analysis of designed nonlinear second-order control systems. As more important the integrators in series, canonical controllable form and Jordan forms are considered. The analysis resumes that due to added controllers systems become stable and insensitive to any disturbance of parameters. Experimental part presents MATLAB simulation of design of possible control systems on the examples of epidemic spread, angular motion of aircraft and submarine depth. The results of simulation confirm the efficiency of offered method of design.

## I. INTRODUCTION

There is a lot of methods of design of robust control which develop with increasing interest and some of them become classical. Commonly all of them are dedicated to defining the ranges of parameters (if uncertainty of parameters takes place) within which the system will function with desirable properties, first of all, will be stable [1,2]. Thus there are many researches which successfully attenuate the uncertain changes of parameters in small (regarding to magnitudes of their own nominal values) ranges. But no one existing method can guarantee the stability of designed control system at arbitrarily large ranges of uncertainly changing parameters of plant. The approach that is offered in the present work relies on the results of catastrophe theory [3,4,5,6,7], uses nonlinear structurally stable functions, and due to bifurcations of equilibrium points in designed nonlinear systems allows to stabilize a dynamic plant with ultimately wide ranges of changing of parameters.

It is known that the catastrophe theory deals with several functions which are characterized by their stable structure. Today there are many classifications of these functions but originally they are discovered as seven basic nonlinearities named as 'catastrophes':

$x^3 + k_1 x$ (fold);
$x^4 + k_2 x^2 + k_1 x$ (cusp);
$x^5 + k_3 x^3 + k_2 x^2 + k_1 x$ (swallowtail);
$x^6 + k_4 x^4 + k_3 x^3 + k_2 x^2 + k_1 x$ (butterfly);
$x_2^3 + x_1^3 + k_1 x_2 x_1 - k_2 x_2 + k_3 x_1$ (hyperbolic umbilic);
$x_2^3 - 3 x_2 x_1^2 + k_1 (x_1^2 + x_2^2) - k_2 x_2 - k_3 x_1$ (elliptic umbilic);
$x_2^2 x_1 + x_1^4 + k_1 x_2^2 + k_2 x_1^2 - k_3 x_2 - k_4 x_1$ (parabolic umbilic).

A part of the catastrophe which does not contain parameters $k_i$ is called as 'germ' of catastrophe. Adding any of them to dynamic system as a controller will give effect shown below. On the example of the catastrophe 'elliptic umbilic' added to dynamical systems we shall see that:
1) new (one or several) equilibrium point appears so there are at least two equilibrium point in new designed system,
2) these equilibrium points are stable but not simultaneous, i.e. if one exists (is stable) then another does not exist (is unstable),
3) stability of the equilibrium points are determined by values or relations of values of parameters of the system,
4) what value(s) or what relation(s) of values of parameters would not be, every time there will be one and only one stable equilibrium point to which the system will attend and thus be stable.

Let us consider the cases of second-order systems (1) and examples (possible applications) of design of control systems (2) of epidemic spread (2-1), aircraft's angular motion (2-2) and submarine depth (2-3).

## II. SECOND ORDER SYSTEMS

**A. Integrators in series.** Let us consider a control plant presented by two integrators connected in series, as shown in Fig.1:

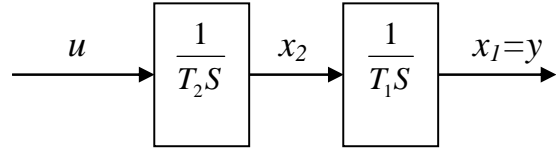

Fig.1. "Integrators in series" structure.

where $T_1$ and $T_2$ are the parameters of integration. The structure of several integrators (more than 2 integrators) is famous of its instability, i.e. no one linear controller can provide the stability to such system and more over with uncertainly changeable parameters [8,9]. The example of two integrators in series allows us to see the advantages of using non-linear catastrophe as controller.

Let us choose a feedback control law as following form:

$$u = -x_2^3 + 3x_2 x_1^2 - k_1\left(x_1^2 + x_2^2\right) + k_2 x_2 + k_3 x_1, \qquad (1)$$

and in order to study stability of the system let us suppose that there is no input signal in the system (equal to zero) [10]. Hence, the system with proposed controller can be presented as:

$$\begin{cases} \dfrac{dx_1}{dt} = \dfrac{1}{T_1} x_2, \\ \dfrac{dx_2}{dt} = \dfrac{1}{T_2}\left(-x_2^3 + 3x_2 x_1^2 - k_1\left(x_1^2 + x_2^2\right) + k_2 x_2 + k_3 x_1\right) \end{cases} \qquad (2)$$

$$y = x_1.$$

The system (2) has following equilibrium points

$$x^1_{1s} = 0, x^1_{2s} = 0; \tag{3}$$

$$x^2_{1s} = \frac{k_3}{k_1}, \quad x^2_{2s} = 0. \tag{4}$$

Stability conditions for equilibrium point (3) obtained via linearization are

$$\begin{cases} -\dfrac{k_2}{T_2} > 0, \\ \dfrac{k_3}{T_1 T_2} < 0. \end{cases} \tag{5}$$

Stability conditions of the equilibrium point (4) are

$$\begin{cases} -\dfrac{3k_3^2 + k_2 k_1^2}{k_1^2 T_2} > 0, \\ \dfrac{k_3}{T_1 T_2} > 0. \end{cases} \tag{6}$$

By comparing the stability conditions given by (5) and (6) we find that the signs of the expressions in the second inequalities are opposite. Also we can see that the signs of expressions in the first inequalities can be opposite due to squares of the parameters $k_1$ and $k_3$ if we properly set their values.

Let us suppose that parameter $T_1$ can be perturbed but remains positive. If we set $k_2$ and $k_3$ both negative and $|k_2| < 3\dfrac{k_3^2}{k_1^2}$ then the value of parameter $T_2$ is irrelevant. It can assume any values both positive and negative (except zero), and the system given by (2) remains stable. If $T_2$ is positive then the system converges to the equilibrium point (3) (becomes stable). Likewise, if $T_2$ is negative then the system converges to the equilibrium point (4) which appears (becomes stable). At this moment the equilibrium point (3) becomes unstable (disappears).

Let us suppose that $T_2$ is positive, or can be perturbed staying positive. So if we can set the $k_2$ and $k_3$ both negative and $|k_2| > 3\dfrac{k_3^2}{k_1^2}$ then it does not matter what value (negative or positive) the parameter $T_1$ would be (except zero), in any case the system (2) will be stable. If $T_1$ is positive then equilibrium point (3) appears (becomes stable) and equilibrium point (4) becomes unstable (disappears) and vice versa, if $T_1$ is negative then equilibrium point (4) appears (become stable) and equilibrium point (3) becomes unstable (disappears).

Results of MatLab simulation for the first and second cases are presented in Fig.2 and 3 respectively. In both cases we see how phase trajectories converge to equilibrium points $(0,0)$ and $\left(\dfrac{k_3}{k_1};0\right)$. In the Fig.2 the phase portrait of the system (2) at constant $k_1=1$, $k_2=-5$, $k_3=-2$, $T_1=100$ and various (perturbed) $T_2$ (from *-4500* to *4500* with step *1000*) with initial condition $x=(-1;0)$ is shown.

In the Fig.3 the phase portrait of the system (2) at constant $k_1=2$, $k_2=-3$, $k_3=-1$, $T_2=1000$ and various (perturbed) $T_1$ (from $-450$ to $450$ with step $100$) with initial condition $x=(-0.25;0)$ is shown.

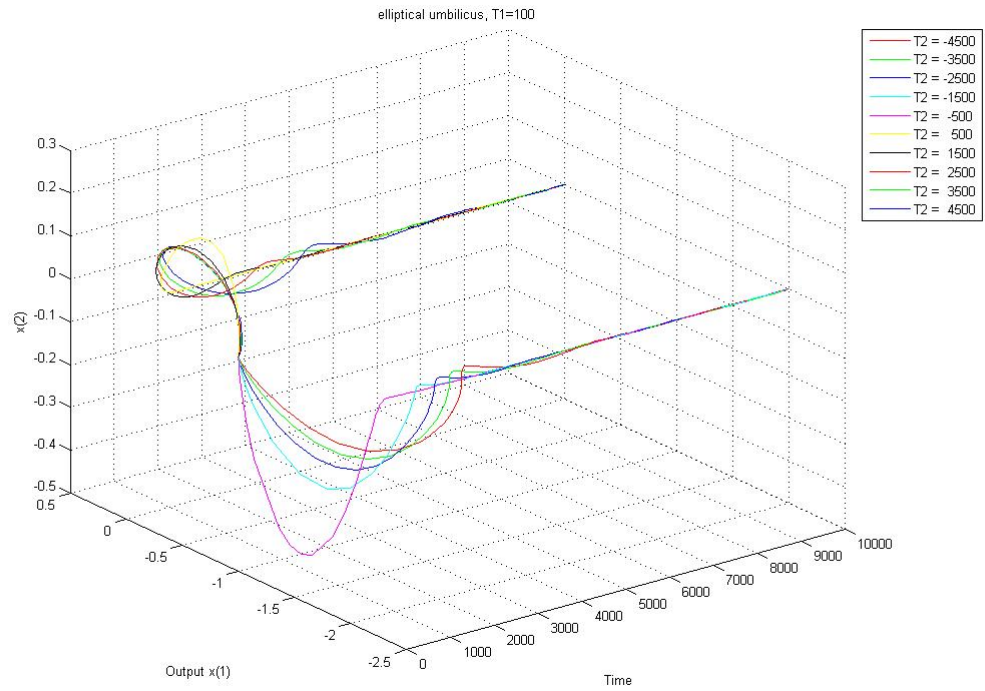

Fig.2. Behavior of output of designed control system in the case of integrators in series at various $T_2$.

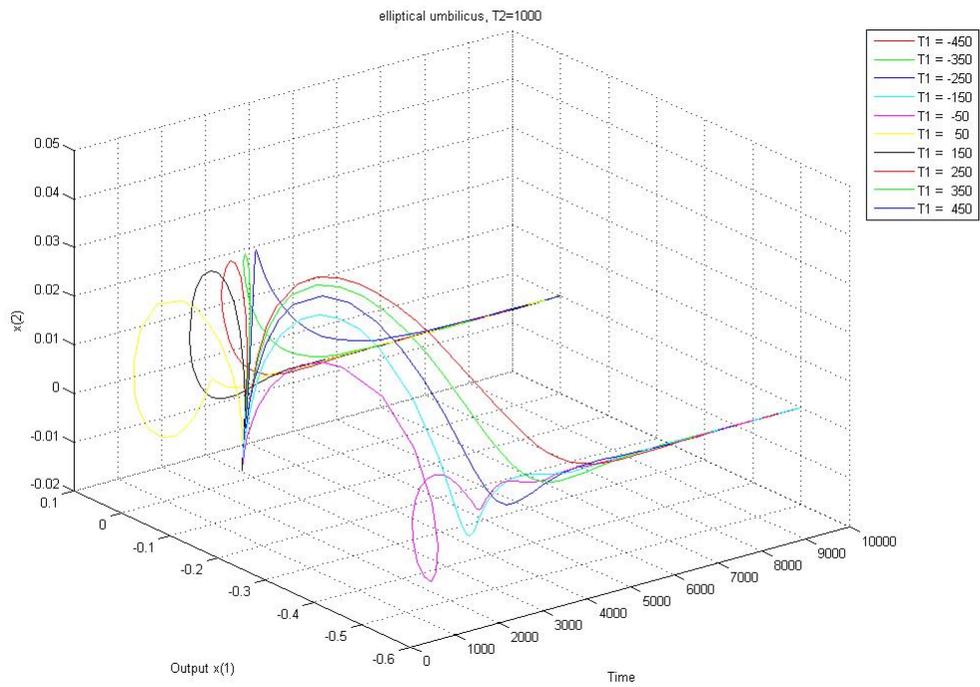

Fig.3. Behavior of output of designed control system in the case of integrators in series at various $T_1$.

Another two forms, canonical controllable form and Jordan form are important because we can reduce any linear matrix of control plant to any of them.

**B. Canonical controllable form (CCF).** This form is important if we would like to affect to the last term of characteristic polynomial $a_n$ which corresponds to general gain of the system.

Let us consider the second order system which is identical to CCF:

$$\begin{cases} \dfrac{dx_1}{dt} = x_2, \\ \dfrac{dx_2}{dt} = -a_2 x_1 - a_1 x_2. \end{cases}$$

$$y = x_1.$$

It is known that the system will be stable if and only if the parameters $a_1$ and $a_2$ are positive. If for example the small perturbation will make the $a_2$ negative then system will become unstable.

Let us set the control law in the form (1). Hence we will obtain the following equations of designed control system.

$$\begin{cases} \dfrac{dx_1}{dt} = x_2, \\ \dfrac{dx_2}{dt} = -a_2 x_1 - a_1 x_2 - x_2^3 + 3x_2 x_1^2 - k_1\left(x_1^2 + x_2^2\right) + k_2 x_2 + k_3 x_1. \end{cases} \quad (7)$$

$$y = x_1.$$

The system (7) has following equilibrium points:

$$x_{1s}^1 = 0, \; x_{2s}^1 = 0; \quad (8)$$

$$x_{1s}^2 = \dfrac{k_3 - a_2}{k_1}, \; x_{2s}^2 = 0; \quad (9)$$

Stability conditions for equilibrium points (8) and (9) respectively are

$$\begin{cases} a_1 - k_2 > 0, \\ a_2 - k_3 > 0. \end{cases} \quad (10)$$

$$\begin{cases} a_1 - k_2 + 3\dfrac{(k_3 - a_2)^2}{k_1^2} > 0, \\ k_3 - a_2 > 0. \end{cases} \quad (11)$$

From inequalities (10) and (11) it is easy to see that here it does not matter what value except zero parameter $a_2$ will be. Similar to above we can resume that system (7) will be stable.

In the Fig.4 the motion of the system (7) at constant control parameters $k_1=4$, $k_2=-4$, $k_3=-6$, constant plant parameter $a_1=1$ and various values of plant parameter $a_2$ which varies from -9.5 to 9.5 with step *1.0*, with initial condition *x=(0.05;0)* is shown.

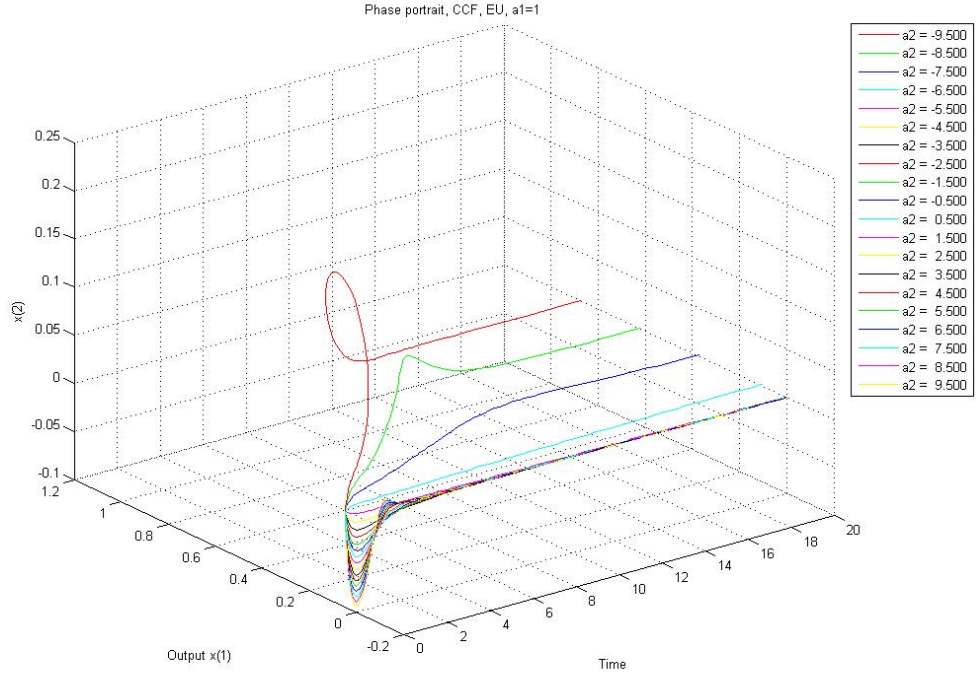

Fig.4. Behavior of output of designed control system in the case of CCF.

**C. Jordan form.** Let us consider the model of second order system which corresponds to Jordan form, i.e. it is a diagonal matrix with eigenvalues as parameters.

$$\begin{cases} \dfrac{dx_1}{dt} = \rho_1 x_1, \\ \dfrac{dx_2}{dt} = \rho_2 x_2. \end{cases} \quad (12)$$

Due to an absence of the relationship between the phase coordinates here we can control each phase coordinate separately and choose the control law in simplified form (as we said without germ).

Let us choose the control law in the simplified form of elliptic umbilic catastrophe without germ and merging the control parameters:

$$u_1 = -k_a x_1^2 + k_b x_1, \; u_2 = -k_a x_2^2 + k_c x_2 \quad (13)$$

Hence, the system (12) with set control (13) is:

$$\begin{cases} \dfrac{dx_1}{dt} = \rho_1 x_1 - k_a x_1^2 + k_b x_1, \\ \dfrac{dx_2}{dt} = \rho_2 x_2 - k_a x_2^2 + k_c x_2. \end{cases} \quad (14)$$

Nonlinear control system (14) has the following equilibrium points:

$$x_{1s}^1 = 0, \; x_{2s}^1 = 0; \tag{15}$$

$$x_{1s}^2 = 0, \; x_{2s}^2 = \frac{\rho_2 + k_c}{k_a}; \tag{16}$$

$$x_{1s}^3 = \frac{\rho_1 + k_b}{k_a}, \; x_{2s}^3 = 0 \tag{17}$$

$$x_{1s}^4 = \frac{\rho_1 + k_b}{k_a}, \; x_{2s}^4 = \frac{\rho_2 + k_c}{k_a}; \tag{18}$$

Stability conditions for the equilibrium point (15) are:

$$\begin{cases} \rho_1 + k_b > 0, \\ \rho_2 + k_c > 0. \end{cases} \tag{19}$$

Stability conditions for the equilibrium point (16) are:

$$\begin{cases} \rho_1 + k_b > 0, \\ \rho_2 + k_c < 0. \end{cases} \tag{20}$$

Stability conditions for the equilibrium point (17) are:

$$\begin{cases} \rho_1 + k_b < 0, \\ \rho_2 + k_c > 0. \end{cases} \tag{21}$$

Stability conditions for the equilibrium point (18) are:

$$\begin{cases} \rho_1 + k_b < 0, \\ \rho_2 + k_c < 0. \end{cases} \tag{22}$$

From inequalities (19-22) it is easy to see that here it does not matter what values except zero parameters $\rho_1$ and $\rho_2$ will be. After any perturbations whatever value except zero this pair would not be, every time there is one and only one of the equilibrium points (15-18) to which the system will attend and at that moment all another equilibrium points will be unstable or will not exist.

In the Fig.5 the motion of the system (14) at constant control parameters $k_a=2$ and $k_b= k_c =5$ and plant various parameters $m_1$ and $m_2$ ($\rho_1$ and $\rho_2$) which vary from *-1250* to *1250* with step *500*, with initial condition *x=(50;50)* is shown.

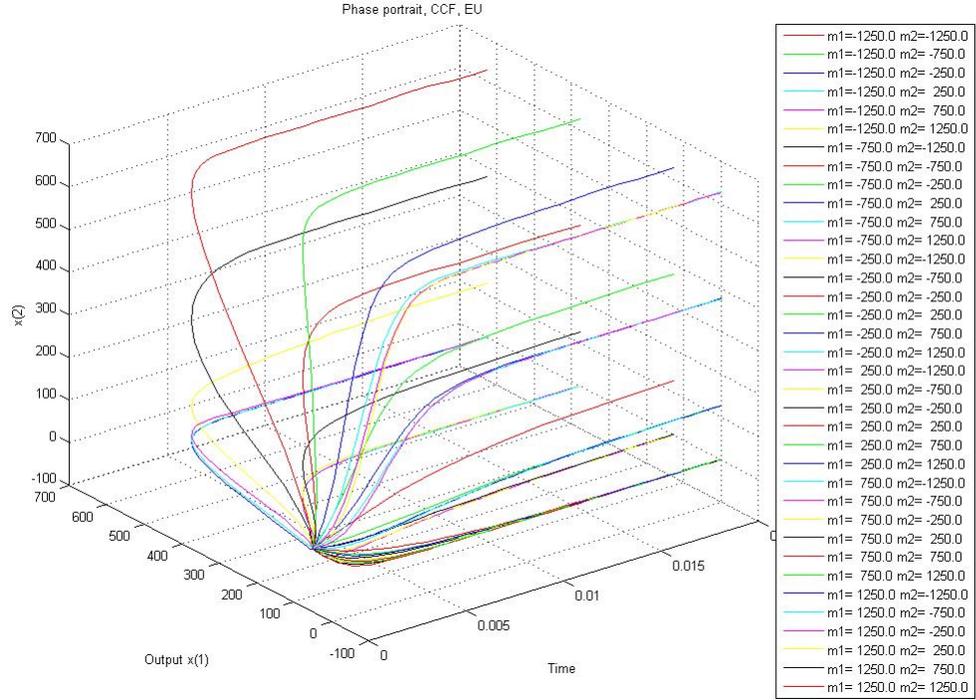

Fig.5. Behavior of output of designed control system in the case of Jordan form.

### III. POSSIBLE APPLICATIONS OF DESIGN OF CONTROL SYSTEMS

**A. Epidemic spread.** The spread of an epidemic disease can be described by a set of differential equations. The population under study is made up of three groups, $x_1$, $x_2$ and $x_3$, such that the group $x_1$ is susceptible to the epidemic disease, group $x_2$ is infected with the disease, and group $x_3$ has been removed from the initial population. The removal of $x_3$ will be due to immunization, death, or isolation from $x_1$. The output of this system is the state variable $x_3$. The plant can be represented by following equations [10]:

$$\begin{cases} \dfrac{dx_1}{dt} = -\alpha x_1 - \beta x_2, \\ \dfrac{dx_2}{dt} = \beta x_1 - \gamma x_2, \\ \dfrac{dx_3}{dt} = \alpha x_1 + \gamma x_2. \end{cases} \quad (23)$$

Let us assume that the population is closed, i.e. the rate at which susceptibles added to the population is equal to 0 and the rate at which new infectives are added to the population is equal 0.

Let us choose the control law in simplified form of catastrophe 'elliptic umbilic' without its germ:

$$u = -k_1\left(x_3^2 + x_2^2\right) + k_2 x_3 + k_3 x_2, \quad B = \begin{pmatrix} 0 \\ 1 \\ 0 \end{pmatrix} \quad (24)$$

Hence, the system (24) with the offered control is:

$$\begin{cases} \dfrac{dx_1}{dt} = -\alpha x_1 - \beta x_2, \\ \dfrac{dx_2}{dt} = \beta x_1 - \gamma x_2 - k_1\left(x_3^2 + x_2^2\right) + k_2 x_3 + k_3 x_2, \\ \dfrac{dx_3}{dt} = \alpha x_1 + \gamma x_2. \end{cases} \quad (25)$$

Hence, the system (25) has two equilibrium points:

$$x_1 = 0,\ x_2 = 0,\ x_3 = 0; \quad (26)$$

$$x_1 = 0,\ x_2 = 0,\ x_3 = \dfrac{k_2}{k_1}. \quad (27)$$

Stability conditions of equilibrium point (26) are:

$$\begin{cases} \alpha + \gamma - k_3 > 0, \\ (\alpha + \gamma - k_3)\left(\alpha(\gamma - k_3) - k_2\gamma + \beta^2\right) - k_2\alpha(\beta - \gamma) > 0, \\ k_2\alpha(\beta - \gamma) > 0. \end{cases} \quad (28)$$

Stability conditions of equilibrium point (27) are:

$$\begin{cases} \alpha + \gamma - k_3 > 0, \\ (\alpha + \gamma - k_3)\left(\alpha(\gamma - k_3) + k_2\gamma + \beta^2\right) + k_2\alpha(\beta - \gamma) > 0, \\ -k_2\alpha(\beta - \gamma) > 0. \end{cases} \quad (29)$$

Here we see the opposite not only for some parameter but for relations of several parameters. To compare the stability with and without offered controller let us see the figures 6 and 7.

Fig.6 presents the output behavior of the system (23) at constant value $\alpha=1$ and various values of parameters $\beta$ and $\gamma$ which vary from 4 to 6 with step 2.

In the Fig.7 the output behavior of the system (25) at constant value $\alpha=1$ and various values of parameters $\beta$ and $\gamma$ which vary from 4 to 6 with step 2 is shown.

As it is proposed the output of the system (25) attends to the values of the equilibrium points depending on the values of parameters $\beta$ and $\gamma$, staying every time stable.

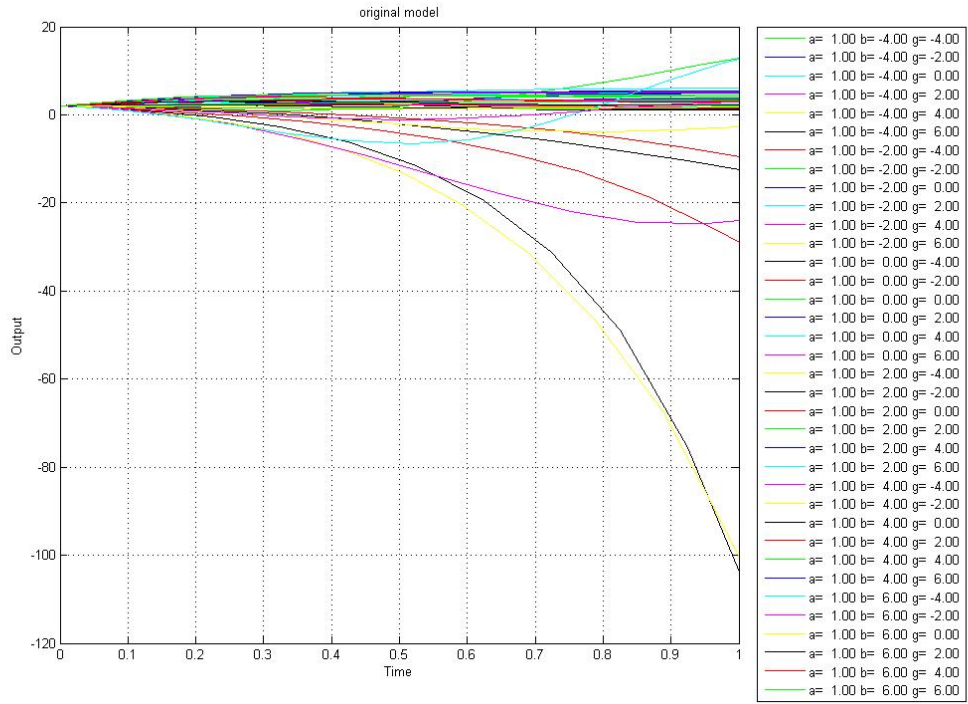

Fig.6. Behavior of rate of number of removals in the epidemic spread at various parameters without controller.

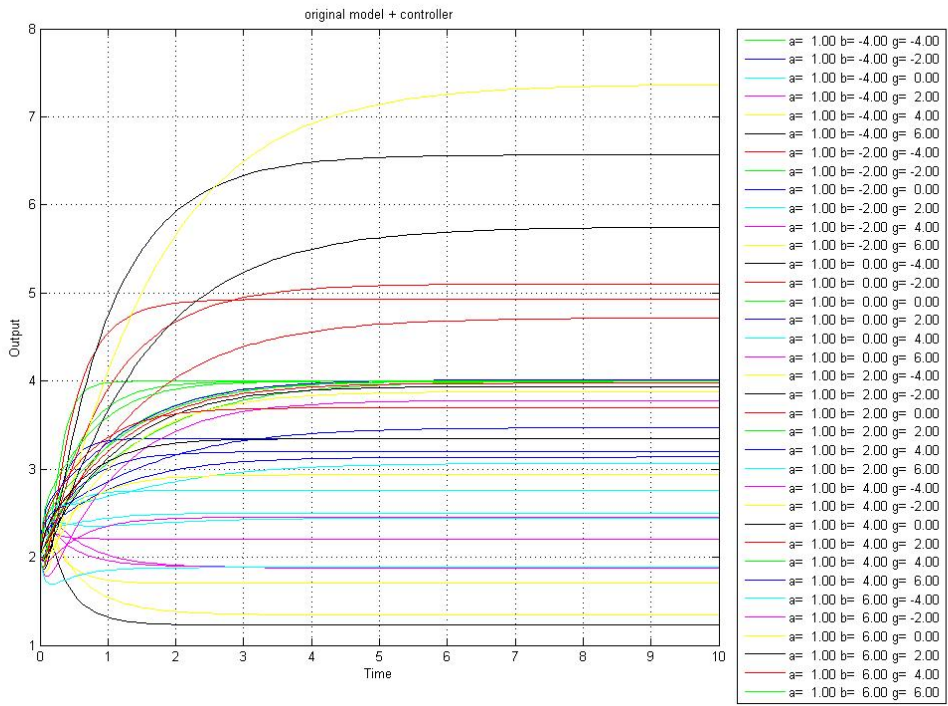

Fig.7. Behavior of rate of number of removals in the epidemic spread at various parameters with controller.

**B. Aircraft's angular motion.** Let us consider the dynamics of aircraft's angular motion. Often it has a quite complicated structure and usually is described by high-order system of nonlinear differential equations [8,9]. But commonly it is possible to isolate a dynamical subsystem which variables and parameters characterize the angles and their relations in attitude of a flight direction as it is shown in the Fig.8 where angle of attack, tilt angle, pitch angle, ground speed, and elevator control signal are denoted as $\alpha$, $\theta$, $\vartheta$, $V$, and $\delta_a$ respectively [9].

Dynamics of aircraft's isolated angular motion is described by the following differential equations:

$$\dot{x} = Ax + Bu,$$
$$y = Cx.$$

where the matrices A, B and C have the following (nominal) values:

$$A = \begin{pmatrix} a_y^\alpha & 0 & -a_y^\alpha \\ a_{m_z}^\alpha & -a_{m_z}^{\omega_z} & -a_{m_z}^\alpha \\ 0 & 1 & 0 \end{pmatrix}, \quad B = \begin{pmatrix} 0 \\ a_{m_z}^{\delta_a} \\ 0 \end{pmatrix}, \quad C = \begin{pmatrix} 0 & 0 & 1 \end{pmatrix}, \tag{30}$$

with the following nominal parameters

$$a_y^\alpha = -2.10 \; [s^{-1}], \; a_{m_z}^\alpha = 29.4 \; [s^{-2}], \; a_{m_z}^{\omega_z} = 2.18 \; [s^{-1}], \; a_{m_z}^{\delta_a} = 60.7 \; [s^{-2}], \; u = \delta_a(t)$$

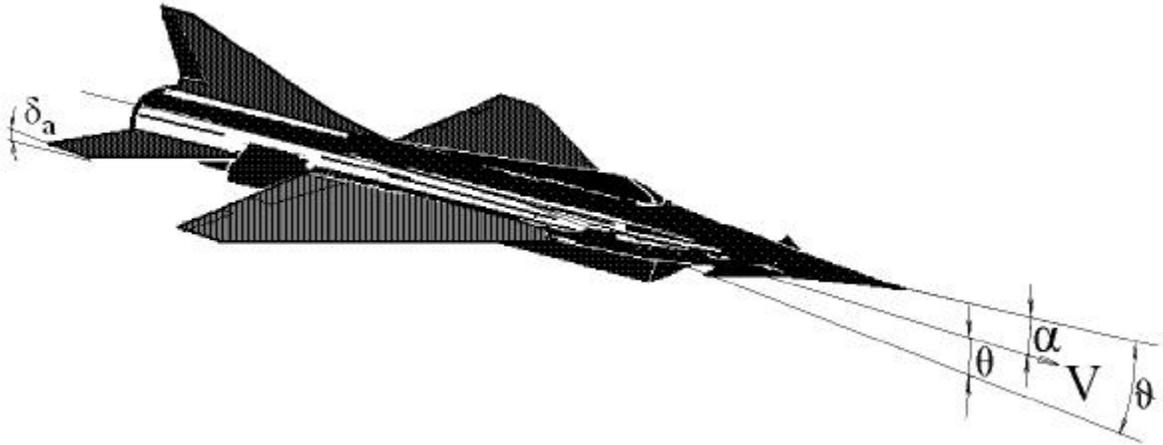

Fig.8. Aircraft's motion characteristics.

If we assume the input $\delta_a(t) = const = 0$ and study the dynamic plant (1) for stability then we see that it is in the stability threshold and not sufficient for engineering practice.

Let us choose controller in the following form:

$$u = -\frac{1}{b_2}\left(k_1(x_3^2 + x_2^2) - k_2 x_3 - k_3 x_2\right). \tag{31}$$

Thus, the system (30) with the added controller (31) will become:

$$\begin{cases} \dfrac{dx_1}{dt} = a_y^\alpha x_1 - a_y^\alpha x_3, \\ \dfrac{dx_2}{dt} = a_{m_z}^\alpha x_1 - a_{m_z}^{\omega_z} x_2 - a_{m_z}^\alpha x_3 - k_1\left(x_3^2 + x_2^2\right) + k_2 x_3 + k_3 x_2, \\ \dfrac{dx_3}{dt} = x_2, \\ y = x_3. \end{cases} \tag{32}$$

New nonlinear control system (32) has two equilibrium points:

$$x_1 = 0,\ x_2 = 0,\ x_3 = 0; \tag{33}$$

$$x_1 = x_3 = \dfrac{k_2}{k_1},\ x_2 = 0. \tag{34}$$

Stability conditions for the equilibrium point (33) are:

$$\begin{cases} a_{m_z}^{\omega_z} - k_3 - a_y^\alpha > 0, \\ \left(a_{m_z}^{\omega_z} - k_3 - a_y^\alpha\right)\left(a_y^\alpha\left(k_3 - a_{m_z}^{\omega_z}\right) - k_2 + a_{m_z}^\alpha\right) - k_2 a_y^\alpha > 0, \\ k_2 a_y^\alpha > 0. \end{cases} \tag{35}$$

Stability conditions for equilibrium point (29) are

$$\begin{cases} a_{m_z}^{\omega_z} + k_3 - a_y^\alpha > 0, \\ \left(a_{m_z}^{\omega_z} + k_3 - a_y^\alpha\right)\left(a_y^\alpha\left(k_3 - a_{m_z}^{\omega_z}\right) + k_2 + a_{m_z}^\alpha\right) + k_2 a_y^\alpha > 0, \\ -k_2 a_y^\alpha > 0. \end{cases} \tag{36}$$

If we draw attention to the last two inequalities in both stability conditions (35) and (36) then we can note the opposite requirements for the sign of the parameter $a_y^\alpha$. Let us assume that the parameter $a_y^\alpha$ satisfies the stability conditions of one of the equilibrium points, i.e. the system converges to that. If after some uncertain perturbation value of the parameter $a_y^\alpha$ is changed such that sign of it becomes opposite then although the current equilibrium point will become unstable or disappear (new value will not satisfy the current stability conditions), another equilibrium point will appear (become stable) because that new value of parameter will automatically satisfy the stability conditions of another equilibrium point. IOW it does not matter for the stability of the system (32) what value except zero the parameter $a_y^\alpha$ would be, in any case the system (32) will be stable.

Let us see the results of MATLAB simulation where one and several parameters varied their values and the system changed phase trajectories but stayed stable.

The Fig.9 shows the behavior of output of the system with added controller at the constant values of the parameters of plant $a_{m_z}^\alpha$=29.4 and $a_{m_z}^{\omega_z}$=2.18, parameter of input $a_{m_z}^{\delta_a}$=60.7 and parameters of control $k_1$=0.1, $k_2$=0.3, $k_3$=0.7 (the values of parameters are chosen arbitrarily) and at the various values of parameter of the plant $a_y^\alpha$ varied from -5.6 to 1.4 with step 0.5 and constant input $\delta_a(t) = const = 1$.

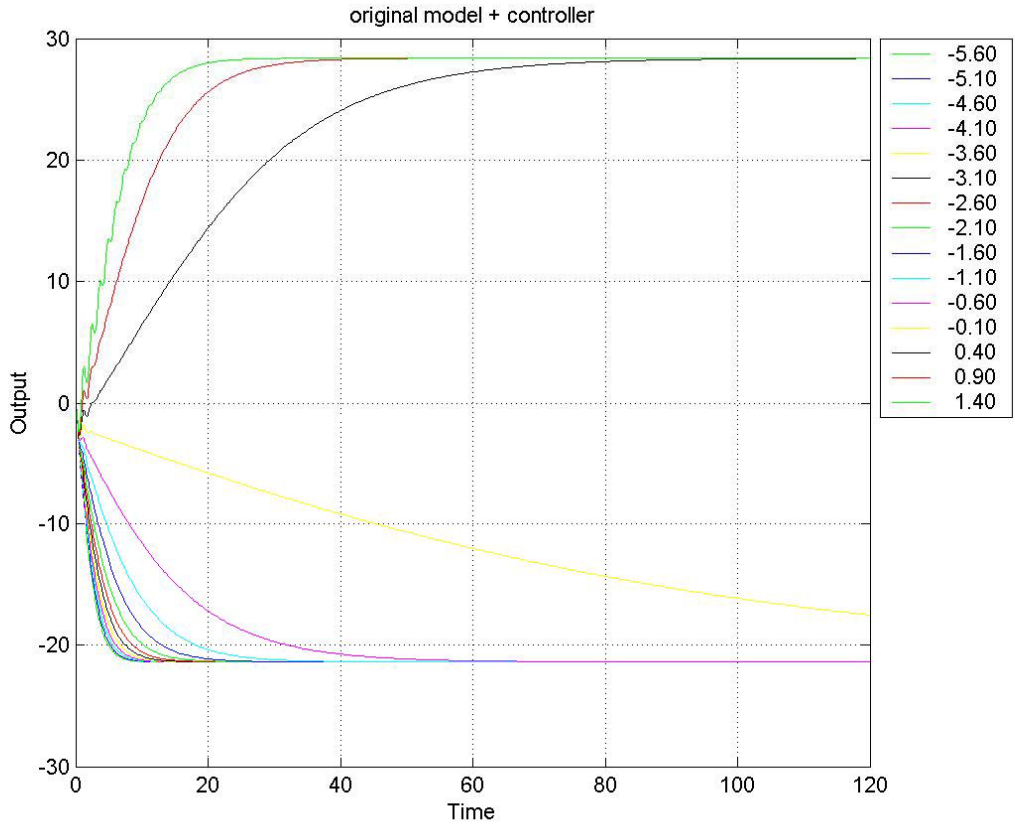

Fig.9. Behavior of output of designed control system of aircraft's angular motion at various $a_y^\alpha$.

In the Fig.10 the behavior of output of the system with added controller at the constant values of the parameters of control $k_1=1$, $k_2=3$, $k_3=7$ (the values of parameters are chosen arbitrarily), parameter of input $a_{m_z}^{\delta_a} = 60.7$ and at the various values of all parameters of the plant $a_y^\alpha$, $a_{m_z}^\alpha$ and $a_{m_z}^{\omega_z}$ varied from -4.1, 9.4 and 0.18 to -0.1, 49.4 and 4.18 (deviations from the nominal values) with steps 1, 10 and 1 respectively and constant input $\delta_a(t) = const = 1$ is shown.

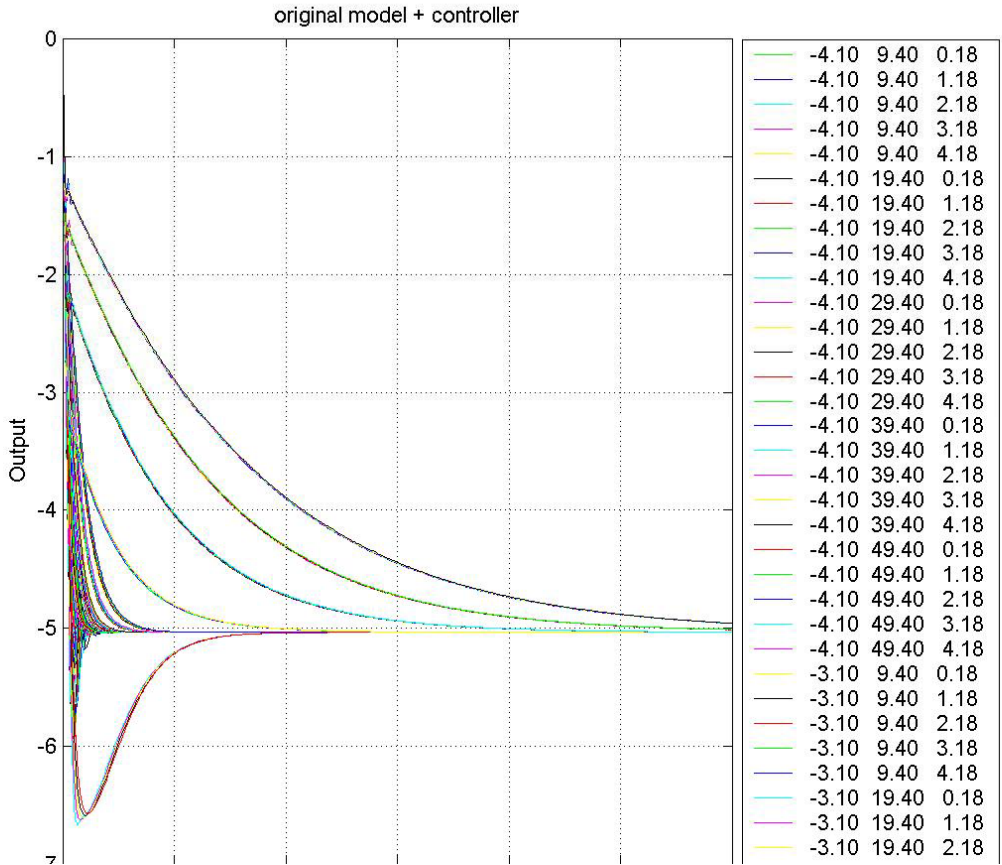

Fig.10. Behavior of output of designed control system of aircraft's angular motion at various $a_y^\alpha$, $a_{m_z}^\alpha$ and $a_{m_z}^{\omega_z}$.

**C. Submarine depth control.** Let us consider dynamics of angular motion of a controlled submarine which is different from the aircraft [10]. This difference results primarily from the moment in the vertical plane due to the buoyancy effect. The important vectors of submarine's motion are shown in the Fig.11.

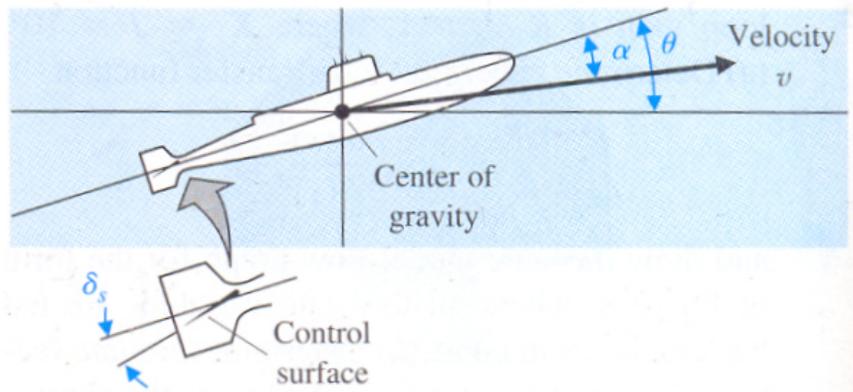

Fig.11. Angles of submarine's depth dynamics.

Let us assume that $\theta$ is a small angle and the velocity $v$ is constant and equal to 25 ft/s. The state variables of the submarine, considering only vertical control, are $x_1 = \theta$, $x_2 = \dfrac{d\theta}{dt}$, $x_3 = \alpha$, where $\alpha$ is the angle of attack and output. Thus the state vector differential equation for this system, when the submarine has an Albacore type hull, is:

$$\dot{x} = Ax + B\delta_s(t), \qquad (37)$$

where

$$A = \begin{pmatrix} 0 & a_{12} & 0 \\ a_{21} & a_{22} & a_{23} \\ 0 & a_{32} & a_{33} \end{pmatrix}, \quad B = \begin{pmatrix} 0 \\ b_2 \\ b_3 \end{pmatrix},$$

parameters of the matrices are equal to:

$a_{12} = 1$, $a_{21} = -0.0071$, $a_{22} = -0.111$, $a_{23} = 0.12$, $a_{32} = 0.07$, $a_{33} = -0.3$,
$b_2 = -0.095$, $b_3 = 0.072$,

and $\delta_s(t)$ is the deflection of the stern plane.

Let us study the behavior of the system (37). In general form it is described as:

$$\begin{cases} \dfrac{dx_1}{dt} = x_2, \\ \dfrac{dx_2}{dt} = a_{21}x_1 + a_{22}x_2 + a_{23}x_3 + b_2\delta_s(t), \\ \dfrac{dx_3}{dt} = a_{32}x_2 + a_{33}x_3 + b_3\delta_s(t). \end{cases} \qquad (38)$$

where input $\delta_s(t)=1$. By turn let us simulate by MATLAB the changing of the value of each parameter deviated from nominal value.

In the Fig.12 the behavior of output of the system (38) at various value of $a_{21}$ (varies from -0.0121 to 0.0009 with step 0.00125) and all left constant parameters with nominal values is presented.

In the Fig.13 the behavior of output of the system (38) at various value of $a_{22}$ (varies from -0.611 to 0.289 with step 0.125) and all left constant parameters with nominal values is presented.

In the Fig.14 the behavior of output of the system (38) at various value of $a_{23}$ (varies from -0.88 to 1.120 with step 0.2) and all left constant parameters with nominal values is presented.

In the Fig.15 the behavior of output of the system (38) at various value of $a_{32}$ (varies from -0.43 to 0.57 with step 0.125) and all left constant parameters with nominal values is presented.

In the Fig.16 the behavior of output of the system (38) at various value of $a_{33}$ (varies from -1.3 to 0.7 to with step 0.25) and all left constant parameters with nominal values is presented.

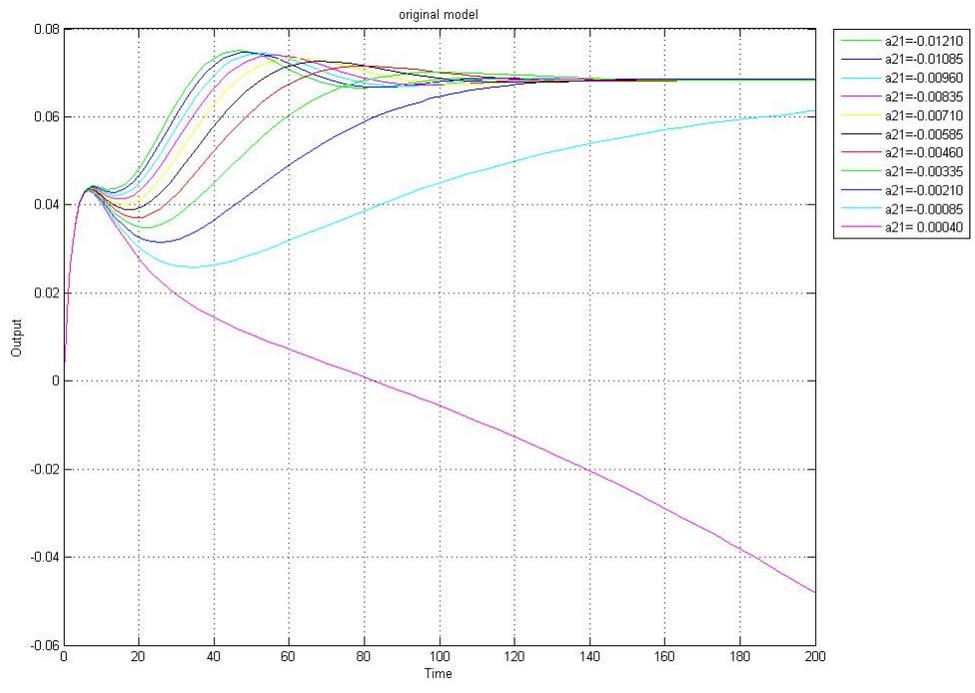

Fig.12. Behavior of output dynamics of submarine's depth at various $a_{21}$.

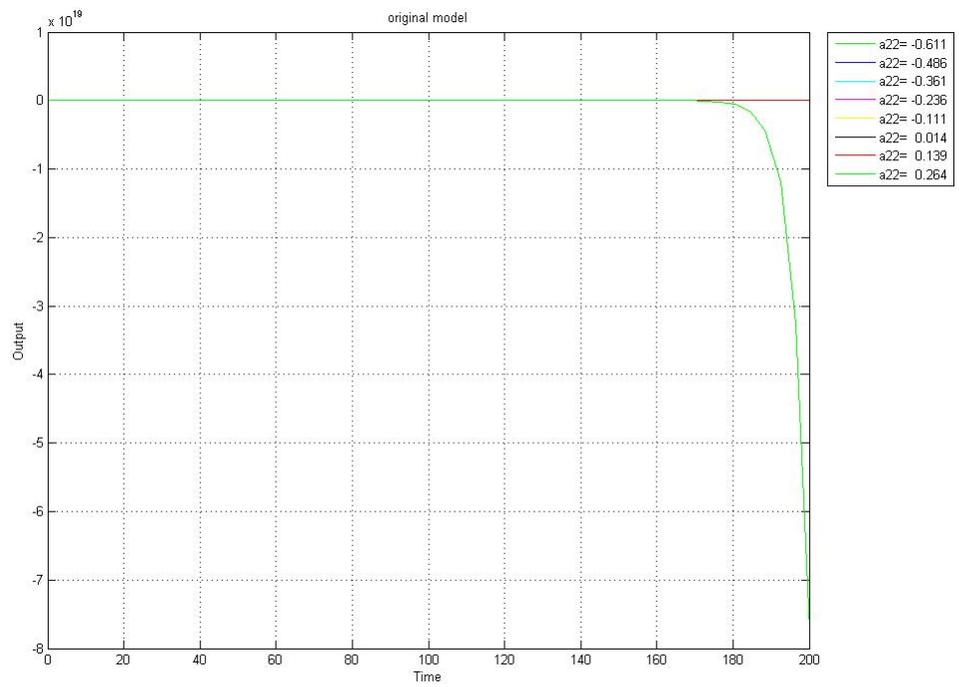

Fig.13. Behavior of output dynamics of submarine's depth at various $a_{22}$.

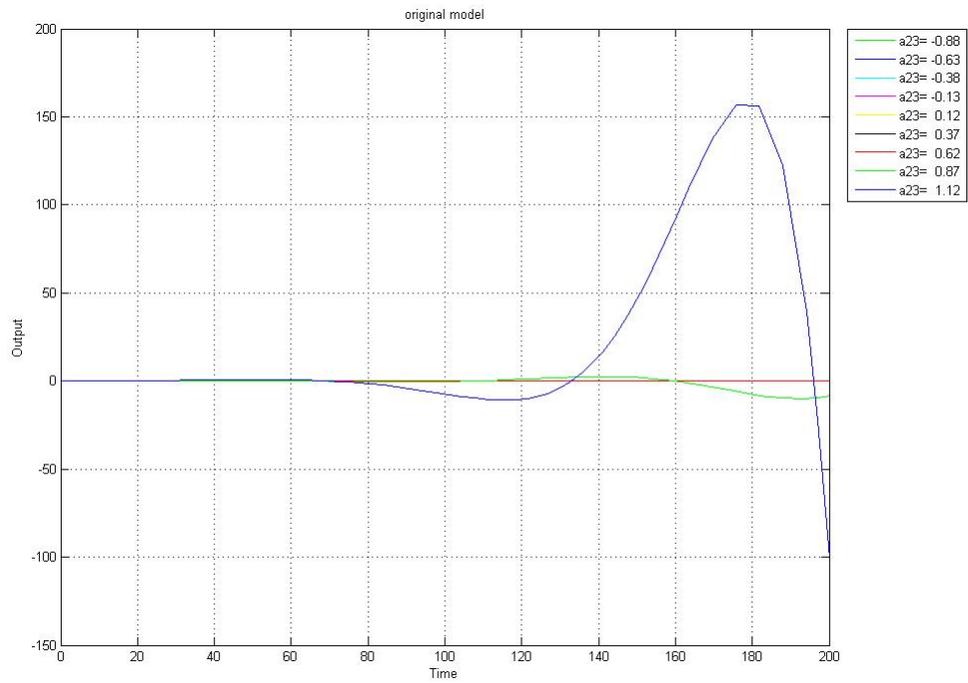

Fig.14. Behavior of output dynamics of submarine's depth at various $a_{23}$.

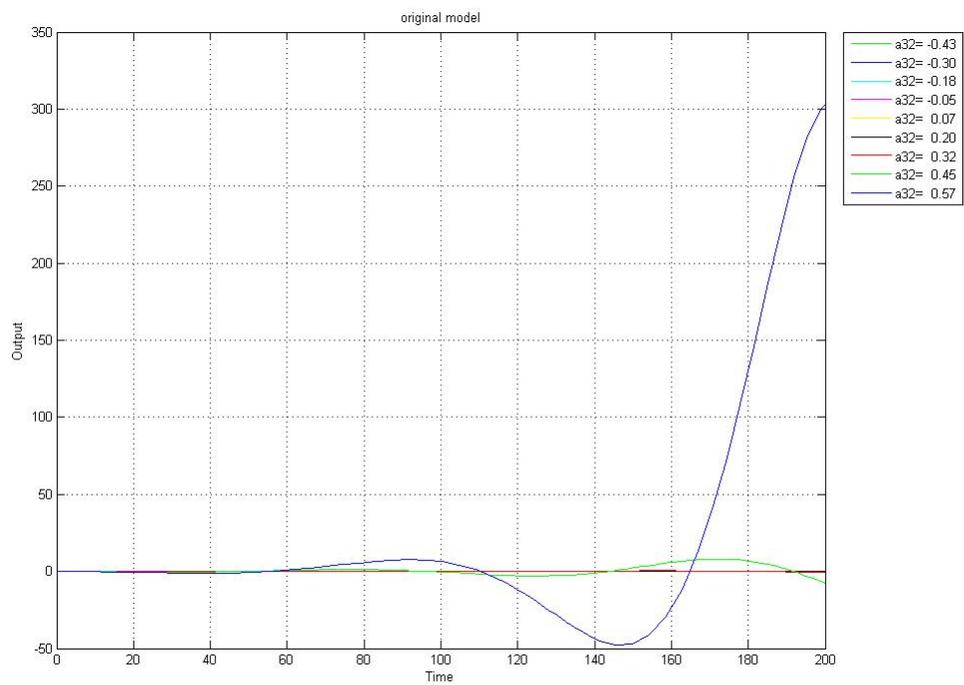

Fig.15. Behavior of output dynamics of submarine's depth at various $a_{32}$.

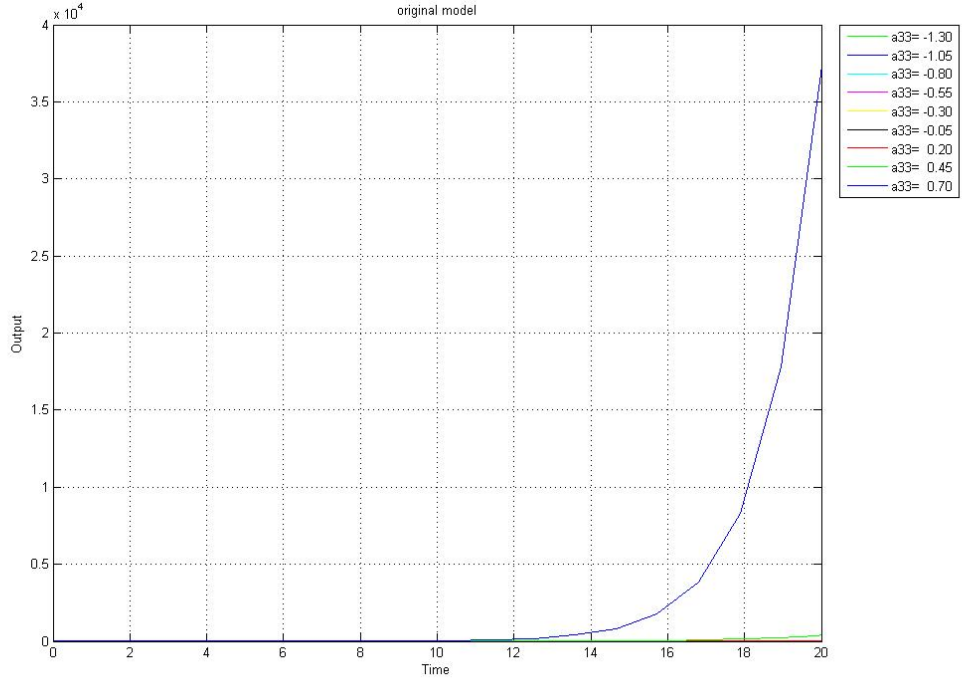

Fig.16. Behavior of output dynamics of submarine's depth at various $a_{33}$.

It is clear that the perturbation of only one parameter makes the system unstable.
Let us set the feedback control law in the following form:

$$u = -k_1(x_3^2 + x_2^2) + k_2 x_3 + k_3 x_2. \tag{39}$$

Hence, designed control system is:

$$\begin{cases} \dfrac{dx_1}{dt} = x_2, \\ \dfrac{dx_2}{dt} = a_{21} x_1 + a_{22} x_2 + a_{23} x_3 + b_2 \delta_S(t), \\ \dfrac{dx_3}{dt} = a_{32} x_2 + a_{33} x_3 + b_3 \delta_S(t) - k_1(x_2^2 + x_3^2) + k_2 x_3 + k_3 x_2. \end{cases} \tag{40}$$

The results of MATLAB simulation of the control system (40) with each changing (disturbed) parameter are presented in the figures 17, 18, 19, 20, and 21.

In the Fig.17 the behavior designed control system (40) at various value of $a_{21}$ (varies from -0.0121 to 0.0009 with step 0.00125) and all left constant parameters with nominal values is presented

In the Fig.18 the behavior of output of the system (40) at various value of $a_{22}$ (varies from -0.611 to 0.289 with step 0.125) and all left constant parameters with nominal values is presented.

In the Fig.19 the behavior of output of the system (40) at various value of $a_{23}$ (varies from -0.88 to 1.120 with step 0.2) and all left constant parameters with nominal values is presented.

In the Fig.20 the behavior of output of the system (40) at various value of $a_{32}$ (varies from -0.43 to 0.57 with step 0.125) and all left constant parameters with nominal values is presented.

In the Fig.21 the behavior of output of the system (40) at various value of $a_{33}$ (varies from -1.3 to 0.7 to with step 0.25) and all left constant parameters with nominal values is presented.

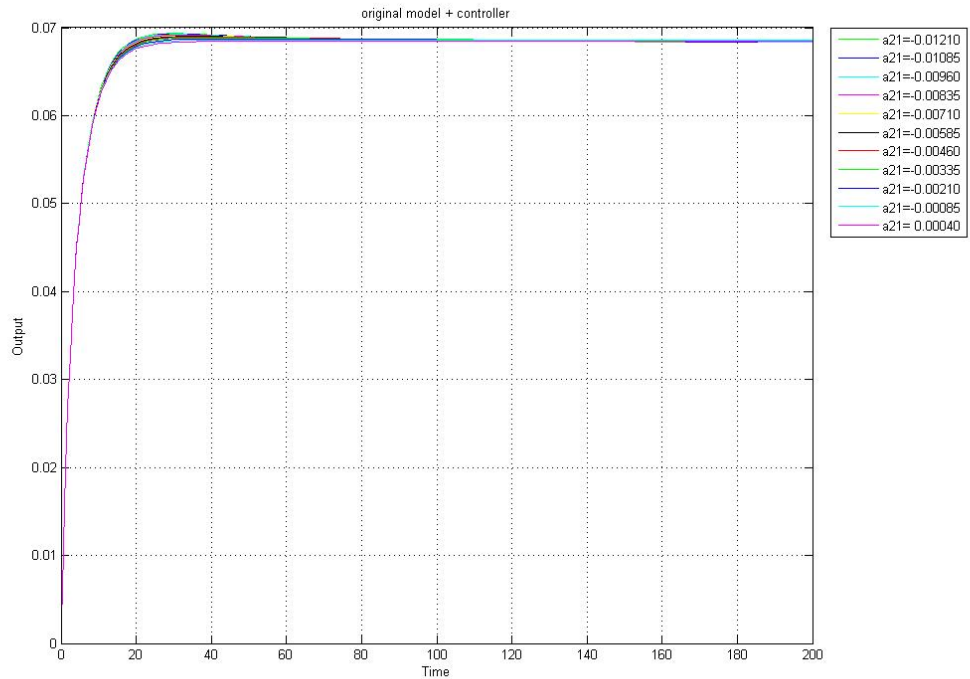

Fig.17 Behavior of output of the submarine depth control system at various $a_{21}$.

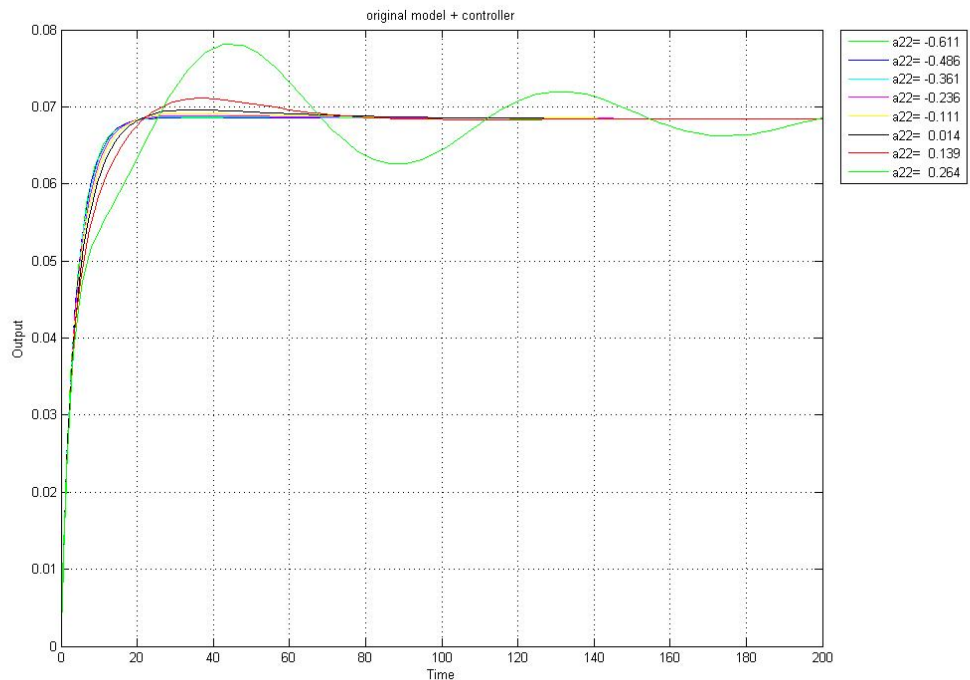

Fig.18. Behavior of output of the submarine depth control system at various $a_{22}$.

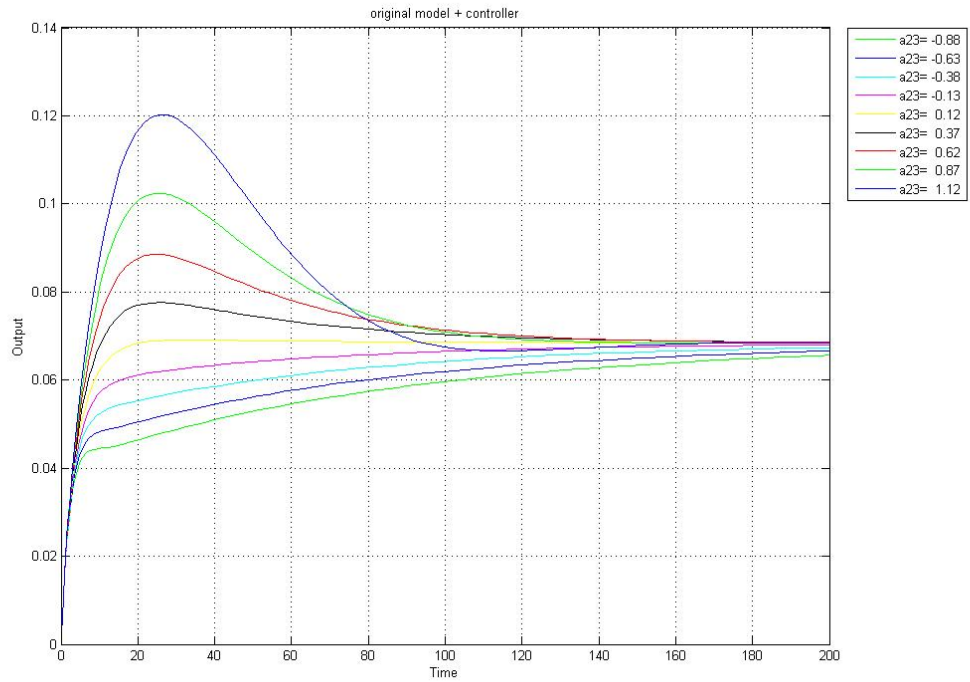

Fig.19. Behavior of output of the submarine depth control system at various $a_{23}$.

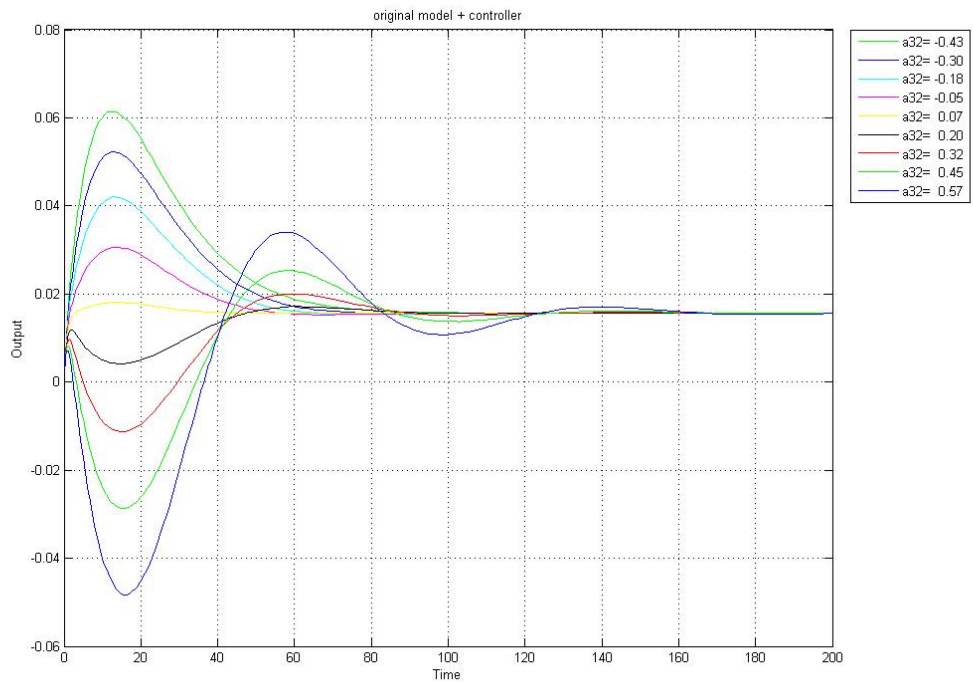

Fig.20. Behavior of output of the submarine depth control system at various $a_{32}$.

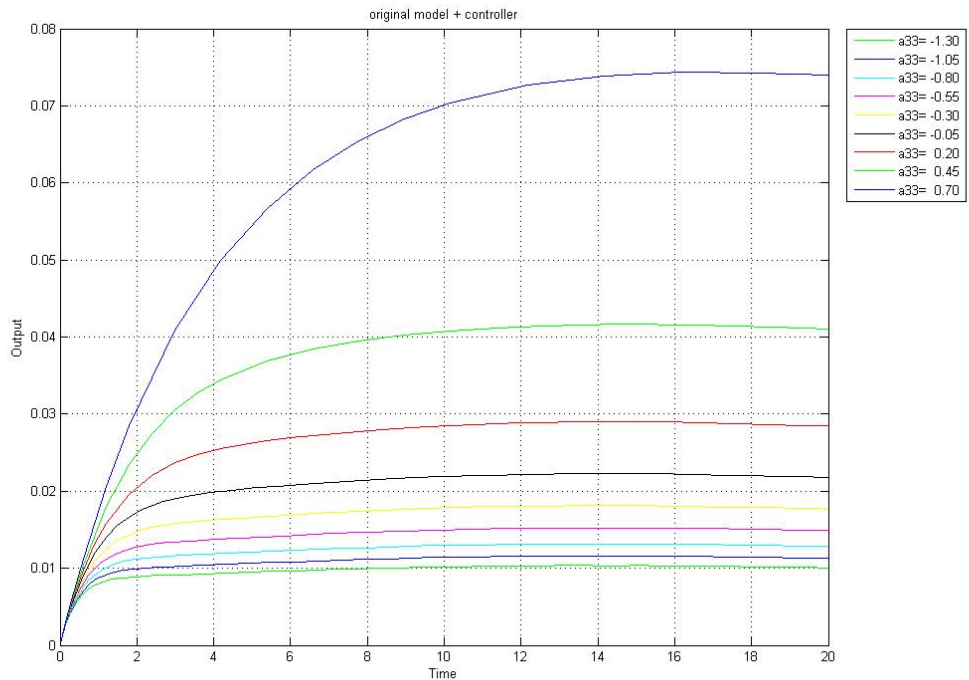

Fig.21. Behavior of output of the submarine depth control system at various $a_{33}$.

IV. Conclusion

Resuming we can conclude that using structurally stable functions from catastrophe theory as controllers give many advantages. The main of them is that the safe ranges of parameters are widened significantly because the designed system stay stable within unbounded ranges of perturbation of parameters even the sign of them changes. The behaviors of designed control systems obtained by MATLAB simulation such that control of epidemic spread, aircraft's angular motion and submarine depth confirm the efficiency of the offered method. The offered approach of design can be applied not only for linear but also for some set or class of nonlinear dynamic plants. For further research and investigation many perspective tasks can occur such that synthesis of control systems with special requirements, design of optimal control, control of chaos, etc.

References

[1] D.-W. Gu, P.Hr. Petkov and M.M. Konstantinov. Robust control design with Matlab. London: Springer-Verlag, 2005.
[2] Boris T. Polyak, Pavel S. Shcherbakov. *Robust stability and control*. Moscow: Nauka, 2002. – 303 pages. (in Russian).
[3] Poston, T. and Stewart, Ian. Catastrophe: Theory and Its Applications. New York: Dover, 1998.
[4] Gilmore, Robert. Catastrophe Theory for Scientists and Engineers. New York: Dover, 1993.
[5] Arnol'd, Vladimir Igorevich. Catastrophe Theory, 3rd ed. Berlin: Springer-Verlag, 1992..
[6] http://en.wikipedia.org/wiki/Catastrophe.
[7] http://en.wikipedia.org/wiki/Catastrophe_theory.
[8] Andrievskii B.R., Fradkov A.L., Izbrannye glavy teorii fvtomaticheskogo upravleniya s primerami na yazyke MATLAB (Selected topics of automatic control theory with examples in the MATLAB language), Petersburg: Nauka, 2000 – 475 p.
[9] Bodner V.A. Aircraft control systems. Moscow: Mashinostroenie, 1973. - 697 p. (in Russian).

[10] Richard C Dorf, Robert H. Bishop. Modern Control Systems, 11/E . Prentice Hall: 2008.
[11] John Doyle, Bruce Francis, Allen Tannenbaum. Feedback control theory. Macmillan Publishing Co., 1990.
[12] Hassan K Khalil. Nonlinear Systems Third Edition Prentice Hall, 2002.
[13] M. Beisenbi, V. Ten. An approach to the increase of a potential of robust stability of control systems. // Theses of the reports of VII International seminar «Stability and fluctuations of nonlinear control systems». Moscow, Institute of problems of control of Russian Academy of Sciences, 2002. - P. 122-123. (in Russian).
[14] Ten V. Design of nonlinear robust controller in a class of structurally stable functions. // Proceedings of World Academy of Science, Engineering and Technology, Volume 5, November 2008, ISSN: 2070370. – p. 523-531.